\documentclass[11pt]{amsart}
\usepackage[T1]{fontenc}
\usepackage{lmodern}
\usepackage[a4paper,margin=29mm]{geometry}
\usepackage{microtype}
\usepackage{amsmath,amssymb,mathtools,mathrsfs}
\numberwithin{equation}{section}
\usepackage{booktabs,longtable,array}
\usepackage{enumitem}
\usepackage{needspace}
\usepackage{placeins}
\usepackage{xcolor}
\usepackage{tikz}
\usetikzlibrary{arrows.meta}
\usepackage{xurl}
\usepackage[colorlinks=true,linkcolor=blue!55!black,citecolor=blue!55!black,urlcolor=blue!65!black]{hyperref}
\usepackage[capitalize,nameinlink]{cleveref}
\setlist{nosep,leftmargin=2.2em}
\allowdisplaybreaks[3]
\newtheorem{theorem}{Theorem}[section]
\newtheorem{lemma}[theorem]{Lemma}
\newtheorem{proposition}[theorem]{Proposition}
\newtheorem{corollary}[theorem]{Corollary}
\theoremstyle{definition}
\newtheorem{definition}[theorem]{Definition}

\theoremstyle{remark}

\crefname{theorem}{theorem}{theorems}
\Crefname{theorem}{Theorem}{Theorems}
\crefname{lemma}{lemma}{lemmas}
\Crefname{lemma}{Lemma}{Lemmas}
\crefname{proposition}{proposition}{propositions}
\Crefname{proposition}{Proposition}{Propositions}
\crefname{corollary}{corollary}{corollaries}
\Crefname{corollary}{Corollary}{Corollaries}
\crefname{definition}{definition}{definitions}
\Crefname{definition}{Definition}{Definitions}
\crefname{example}{example}{examples}
\Crefname{example}{Example}{Examples}
\crefname{remark}{remark}{remarks}
\Crefname{remark}{Remark}{Remarks}
\crefname{appendixlemma}{lemma}{lemmas}
\Crefname{appendixlemma}{Lemma}{Lemmas}
\newcommand{\abs}[1]{\lvert#1\rvert}
\def\dx{\,\mathrm{d}x}

\def\pa{\partial}
\newcommand{\manuscripttitle}{First Order Convergence of Representative Volume Element with Mixed Boundary Condition}
\title[RVE with Mixed Boundary Condition]{\manuscripttitle}
\hypersetup{
  pdftitle={\manuscripttitle},
  pdfauthor={Pingbing Ming; Hao Yu},
  pdfsubject={Periodic homogenization and mixed Dirichlet--Neumann cell problems},
  pdfkeywords={periodic homogenization, effective matrix, mixed boundary conditions, nonsymmetric elliptic equations, measurable coefficients}}

\author[P. Ming]{Pingbing Ming}
\author[H. Yu]{Hao Yu}
\address{SKLMS, Institute of Computational Mathematics and Scientific/Engineering Computing, Academy of Mathematics and Systems Science, Chinese Academy of Sciences,
Beijing 100190, China; School of Mathematical Sciences, University of Chinese Academy of Sciences, Beijing 100049, China}
\email[Pingbing Ming]{mpb@lsec.cc.ac.cn}
\email[Hao Yu]{yuhao@amss.ac.cn}

\begin{document}
\begin{abstract}
We study the cell problem with mixed Dirichlet--Neumann boundary condition to reconstruct the effective matrix for two-dimensional scalar divergence-form equations with bounded measurable, and possibly nonsymmetric coefficients. We prove uniform invertibility of the mean-gradient matrix. Moreover, we derive the first-order error estimate when the medium is periodic.
\end{abstract}

\maketitle

\section{Introduction}
\label{sec:introduction}

A central problem in many multiscale methods is to find the macroscopical properties of the solution. This is achieved by solving a lot of local problems posed on computational domains, which are called cells in the literature. There are many methods towards this objective, such as numerical upscaling~\cite{Durlofsky:1991,PecullanGibianskyTorquato1999}, heterogeneous multiscale method~\cite{EE:2003}, among many others. A typical example is the representative volume element (RVE)~\cite{Durlofsky:1991}. Key ingredients of such methods, particularly for RVE, are the boundary conditions and the size of the cell. The boundary conditions are crucial to the accuracy and the efficiency of the method. Though they seem artificial at first glance, they are required only because the cells are truncated and localized.

We confine ourselves to the following problem, with \(a^\epsilon(x)=a(x/\epsilon)\) and the periodic coefficient \(a\) specified below:
\[
-\nabla\cdot(a^\epsilon(x)\nabla u^{\epsilon}(x))=f(x),\quad x\in\Omega\subset\mathbb{R}^2,\qquad u^{\epsilon}(x)=0,\quad x\in\partial\Omega,
\]
where $\epsilon$ is a small parameter that signifies the multiscale nature of the coefficient family \(a^\epsilon\). We are interested in the accuracy of the method with which the effective matrix is obtained. For the periodic coefficient, the effective matrix \(\mathcal A\) is given explicitly by the corrector formula~\eqref{eq:H}.

In practice, the effective matrix is computed by solving an array of local problems over a finite sample block. The resulting effective matrix depends on the sample size and on the imposed boundary conditions. For a sampling cell of size \(\delta\), boundary interactions with the oscillations typically produce a resonance error of order \(\epsilon/\delta\); see~\cite{EMingZhang2005,DuMing2010,Gloria2011}. The same finite size effect underlies representative volume approximations, where random media introduce the additional distinction between fluctuations and systematic bias~\cite{BourgeatPiatnitski2004,Fischer2019,SchneiderJosienOtto2022}. We focus on a periodic medium and seek an estimate valid for every microscopic phase.

The study on the effect of the boundary conditions on such local problems has also been reported in the literature. \textsc{Durlofsky}~\cite{Durlofsky:1991} has studied the extraction of the effective permeability of a porous medium by solving Darcy's equation over a representative elementary volume and has considered the influence of different boundary conditions. Particularly, he investigated the numerical performance of different formulations of mixed Dirichlet--Neumann boundary conditions. \textsc{Pecullan et al}~\cite{PecullanGibianskyTorquato1999} considered the effective elastic properties for the two-dimensional composites. \textsc{Yue and E} systematically studied the influence of the boundary conditions in~\cite{YueE2007}. They pointed out that several formulations for the RVE with mixed Dirichlet--Neumann boundary conditions are incorrect. We refer to~\cite{RenardAbabou2022} for a state-of-the-art review on this topic. In the theoretical part, \textsc{Wu, Efendiev and Hou}~\cite{WuEfendievHou2002} have derived the convergence rate for RVE with Dirichlet boundary condition. The same result has been obtained in~\cite{EMingZhang2005} in the framework of the heterogeneous multiscale method. \textsc{Du and Ming} unified the proof for the convergence rate for RVE with different boundary conditions, which covers the Dirichlet boundary condition, the Neumann boundary conditions and the generalized periodic boundary conditions. Unfortunately, their method does not work for the RVE with mixed boundary condition. The purpose of this work is to prove the convergence rate for RVE with mixed boundary conditions for two-dimensional periodic media.

The local problem studied here imposes a pressure drop in one coordinate direction and zero conormal flux in the transverse direction. The average fluxes and gradients of the two coordinate loadings form matrices \(F\) and \(G\), respectively, and the approximated effective matrix is reconstructed as \(\mathcal{A}_{\eta}{:}=FG^{-1}\). For a physical cell of side length \(\delta\), the dimensionless parameter is \(\eta=\epsilon/\delta\).

Our main result establishes both stability and consistency for two-dimensional elliptic problems in divergence form. For every bounded measurable coefficient satisfying the ellipticity assumptions in \eqref{eq:E}, and every \(0<\eta\le1\), we prove
\(\sigma_{\min}(G_\eta)\ge c(\lambda,\Lambda)>0\). Moreover, if the coefficient is periodic, we obtain
\(\|\mathcal A_\eta-\mathcal A\|_2\le C(\lambda,\Lambda)\eta\).
On a physical cell this gives the rate \(\mathcal{O}(\epsilon/\delta)\). For every nondegenerate ellipticity range \(\lambda<\Lambda\), an explicit smooth coefficient attains a matching first-order error along a sequence of noninteger cell sizes. Thus the rate is optimal when arbitrary cell sizes are allowed, although commensurate cells may be exact.

The remaining part of the paper is as follows. In~\Cref{sec:problem}, we state the main results. In~\Cref{sec:prelim,sec:invertibility,sec:stream} the well-posedness, uniform invertibility, and the planar dual problem have been established. In~\Cref{sec:boundary-functionals}, we derive the boundary representation, and \cref{sec:completion} proves the estimate and sharpness of the error estimate.
\section{Mixed cell problem and the main result}
\label{sec:problem}
We first fix certain notations. For any bounded Lipschitz domain $\Omega$, the standard Sobolev space $H^s(\Omega)$ will be used~\cite{AdamsFournier2003} with $s\in\mathbb{R}$. For $1/2<s<3/2$, the trace map $\operatorname{Tr}:H^s(\Omega)\to H^{s-1/2}(\partial\Omega)$ is bounded. We will use $\abs{\cdot}$ to denote the absolute value of a scalar quantity, the Euclidean norm of a vector and the volume of a set $K$. 

For \(0<\lambda\le\Lambda<\infty\), define the pointwise matrix class
\begin{equation}
\mathbb M(\lambda,\Lambda):=
\bigl\{B\in\mathbb R^{2\times2}\mid\,
\xi\cdot B\xi\ge\lambda|\xi|^2,\,\xi\cdot B\xi\ge\Lambda^{-1}|B\xi|^2\;
\text{ for all\;}\xi\in\mathbb R^2\bigr\}.
\label{eq:E}
\end{equation}
For a measurable set \(D\subset\mathbb R^2\), let
\(\mathcal M(\lambda,\Lambda;D)\) be the class of measurable fields
\(b:D\to\mathbb R^{2\times2}\) such that
\(b(x)\in\mathbb M(\lambda,\Lambda)\) for almost every \(x\in D\).
We abbreviate by \(\mathcal M(\lambda,\Lambda)\) the subclass of
\(\mathcal M(\lambda,\Lambda;\mathbb R^2)\) consisting of
\(\mathbb Z^2\)-periodic fields. Throughout, we assume
\(a\in\mathcal M(\lambda,\Lambda)\).

Let \(Q=(0,1)^2\), and denote its four open edges by \(\Gamma_L,\Gamma_R,\Gamma_B,\Gamma_T\), corresponding to the left, right, bottom, and top edges, respectively. Define
\[
\begin{aligned}
V_1{:}&=\{w\in H^1(Q)\,\mid\,\operatorname{Tr}w=0\text{ on }\Gamma_L\cup\Gamma_R\},\\
V_2{:}&=\{w\in H^1(Q)\,\mid\,\operatorname{Tr}w=0\text{ on }\Gamma_B\cup\Gamma_T\}.
\end{aligned}
\]

Fix \(0<\eta\le1\) and \(\vartheta\in\mathbb R^2\), and set \(a_\eta(x):=a(\vartheta+x/\eta)\) for \(x\in Q\). For \(j=1,2\), find \(u_j\) with \(u_j-x_j\in V_j\) such that
\begin{equation}
\int_Q a_\eta\nabla u_j\cdot\nabla\phi\dx=0
\quad\text{for every }\phi\in V_j.
\label{eq:W}
\end{equation}
Denote the flux by \(q_j=a_\eta\nabla u_j\), and the energy by \(\mathscr E_j=\int_Qq_j(x)\cdot\nabla u_j(x)\dx\). 

The cell averages are defined by
\begin{equation}
\begin{aligned}
G_\eta&=\left[\;\int_Q\nabla u_1(x)\dx
\;\middle|\;\int_Q\nabla u_2(x)\dx\;\right],\\
F_\eta&=\left[\;\int_Qa_\eta(x)\nabla u_1(x)\dx
\;\middle|\;\int_Qa_\eta(x)\nabla u_2(x)\dx\;\right].
\end{aligned}
\label{eq:FG}
\end{equation}

The matrices collect the two average gradients and fluxes as columns. Whenever \(G_\eta\) is invertible, define
\(\mathcal A_\eta:=F_\eta G_\eta^{-1}\). The main theorem below shows that this reconstruction is well defined uniformly in \(\eta\) and \(\vartheta\).

The target matrix is defined independently through the zero-mean periodic correctors \(\chi_j\in H^1_{\rm per}(Q)\) with \(\int_Q\chi_j=0\):
\begin{equation}
\begin{aligned}
\int_Q a(y)(e_j+\nabla\chi_j(y))\cdot\nabla\phi(y)\,\mathrm{d}y
&=0,\qquad\text{for all\quad}\phi\in H^1_{\rm per}(Q),\\
\mathcal A e_j&=\int_Qa(y)(e_j+\nabla\chi_j(y))\,\mathrm{d}y.
\end{aligned}
\label{eq:H}
\end{equation}

Define
\[
c_{\ast}{:}=\frac14\exp\left(-1-\frac{\pi\Lambda}{2\lambda}\right).
\]

We state the main result as follows.
\begin{theorem}\label{thm:main}
For every 
\(\eta\in(0,1]\) and \(\vartheta\in\mathbb R^2\),
\[
\sigma_{\min}(G_\eta)\ge c_*.
\]

Moreover, there exists a constant \(C\), depending only on \(\lambda\) and \(\Lambda\), such that, for every such \(\eta\) and \(\vartheta\),
\[
\|\mathcal A_\eta-\mathcal A\|_2\le C\eta.
\]
\end{theorem}

We transform the above error estimate to the physical cell:
\begin{corollary}
\label{cor:physical_rate}
Let $a^{\epsilon}(x)=a(x/\epsilon)$ with $\{a_{ij}\}_{i,j=1}^2$ $1-$periodic functions. Denote \(I_\delta=x_*+\delta Q\), where \(0<\epsilon\le\delta\), and set \(V_j(I_\delta):=\{\phi\in H^1(I_\delta):\phi(x_*+\delta\,\cdot)\in V_j\}\). For each \(j=1,2\), let \(v_j^\epsilon-x_j\in V_j(I_\delta)\) solve \(\int_{I_\delta}a^\epsilon\nabla v_j^\epsilon\cdot\nabla\phi\,dx=0\) for every \(\phi\in V_j(I_\delta)\). With \(\langle g\rangle_{I_\delta}:=|I_\delta|^{-1}\int_{I_\delta}g\,dx\), set \(G_H:=[\langle\nabla v_1^\epsilon\rangle_{I_\delta}\mid\langle\nabla v_2^\epsilon\rangle_{I_\delta}]\) and \(F_H:=[\langle a^\epsilon\nabla v_1^\epsilon\rangle_{I_\delta}\mid\langle a^\epsilon\nabla v_2^\epsilon\rangle_{I_\delta}]\). %
Define \(\mathcal A_H:=F_HG_H^{-1}\). Then
\begin{equation}
\|\mathcal A_H-\mathcal A\|_2
\le C(\lambda,\Lambda)\frac{\epsilon}{\delta},
\qquad 0<\epsilon\le\delta.
\label{eq:PHYS}
\end{equation}
\end{corollary}

Here \(\delta/\epsilon\) is the number of microscopic periods across the cell.
The order in \cref{eq:PHYS} cannot be improved uniformly when arbitrary, noninteger cell sizes are admitted.

\begin{proposition}\label{prop:sharpness}
For every \(0<\lambda<\Lambda\), set \(m=(\lambda^{-1}+\Lambda^{-1})/2\) and \(r=(\lambda^{-1}-\Lambda^{-1})/4>0\). There exist a fixed smooth coefficient \(a\in\mathcal M(\lambda,\Lambda)\) and a sequence \(\eta_N\to 0\) such that
\begin{equation}
\|\mathcal A_{\eta_N}-\mathcal A\|_2
=\frac{r\eta_N}{m(m\pi+r\eta_N)}.
\label{eq:sharpness-exact}
\end{equation}
\end{proposition}

\section{Main Technical Results}\label{sec:main}
\subsection{Matrix preliminaries, weak solutions, and energy bounds}
\label{sec:prelim}
We list certain properties of the transpose and rotational duality for elements in \(\mathbb M(\lambda,\Lambda)\).
\begin{lemma}\label{lem:matrix_tools}
For \(b\in\mathbb M(\lambda,\Lambda)\), the following properties hold:
\begin{enumerate}[label=\textup{(\roman*)}]
\item The matrix \(b\) is invertible, \(\|b\|_2\le\Lambda\), and \(\|b^{-1}\|_2\le1/\lambda\).
\item The transpose satisfies \(b^{\top}\in\mathbb M(\lambda,\Lambda)\).
\item With \(J:=\begin{pmatrix}0&-1\\1&0\end{pmatrix}\), the dual matrix
\(b^\#{:}=-Jb^{-1}J\) belongs to \(\mathbb M(1/\Lambda,1/\lambda)\).
\end{enumerate}
\end{lemma}

The assertions follow directly from~\eqref{eq:E}. The transformation
\(b^\#=J^{-1}b^{-1}J\) is the standard planar dual
conductivity~\cite[(3) and the subsequent identity]{BrianeManceau2008};
for these ellipticity conditions and their invariance under transposition, see also~\cite[Proposition 2]{Muratartar1997}.

For the remainder of this section, fix \(b\in\mathcal M(\lambda,\Lambda;Q)\), without assuming periodicity. Write \(u_j,q_j,\mathscr E_j,F,G\) for the mixed solutions, fluxes, energies, and average matrices defined in \Cref{sec:problem} with \(a_\eta\) replaced by \(b\), omitting the subscript \(\eta\) on \(F,G\).
\begin{lemma} %
\label{lem:rectangle_traces}
The trace map \(\operatorname{Tr}:H^1(Q)\to L^2(\partial Q)\) is bounded. Moreover, for almost every coordinate line, the corresponding slice belongs to \(H^1(0,1)\), its one-dimensional weak derivative is the corresponding partial derivative, and its endpoint values agree with the two-dimensional trace.
\end{lemma}

\begin{proof}
This is the specialization to the square of Evans--Gariepy~\cite[Theorems~4.6, 4.7, 4.20, and~4.21]{EvansGariepy2015}. The extension theorem reduces the slicing assertion to the whole-space characterization of \(H^1(\mathbb R^2)\)-functions by representatives that are absolutely continuous on almost every line parallel to a coordinate axis, while approximation identifies the one-dimensional endpoint values with the Sobolev trace, which completes the proof.
\end{proof}

We shall also use the coordinate-slice form of the Friedrichs inequality. We omit the proof for simplicity.
\begin{lemma}
\label{lem:mixed_space_poincare}
If \(f\in H^1(0,1)\) and \(f(0)=0\) or \(f(1)=0\), then
\(\|f\|_{L^2(0,1)}^2\le\frac12\|f'\|_{L^2(0,1)}^2\).
Moreover, for \(j=1,2\) and \(w\in V_j\),
\[
\|w\|_{L^2(Q)}^2
\le\frac12\|\partial_jw\|_{L^2(Q)}^2
\le\frac12\|\nabla w\|_{L^2(Q)}^2.
\]
\end{lemma}

Consequently, for \(j=1,2\), \(V_j\) is a Hilbert space under the inner product
\[
(w,\phi)_{V_j}:=\int_Q\nabla w(x)\cdot\nabla\phi(x)\dx.
\]

The following lemma yields certain refined properties for the solution of the cell problems.
\begin{lemma}\label{lem:wellposed_energy}
The two mixed problems in~\eqref{eq:W} have unique solutions $u_j$ satisfying
\begin{gather}
\int_Q\partial_ju_j(x)\dx=1,\label{eq:B1}
\\
\lambda\le\mathscr E_j\le\Lambda,\qquad
\|\nabla u_j\|_{L^2(Q)}^2\le\Lambda/\lambda,
\qquad \|q_j\|_{L^2(Q)}\le\Lambda.
\label{eq:B2}
\end{gather}
Moreover, \(0\le u_j\le1\) almost everywhere in \(Q\) for \(j=1,2\).
\end{lemma}

\begin{proof}
By Lemma~\ref{lem:mixed_space_poincare}, the existence and uniqueness of $u_j$ are standard, and we omit the proof for simplicity.

An integration by parts gives~\eqref{eq:B1}. 

By construction, \(u_j-x_j=w\in V_j\). Thus \(u_j-x_j\) is an admissible test function in~\eqref{eq:W}. A direct calculation gives \(\mathscr E_j=\int_Qe_j\cdot q_j(x)\dx\), which together with the Cauchy--Schwarz inequality leads to \(\mathscr E_j\le\|q_j\|_2\). By ellipticity, we obtain
\(\|q_j\|_2^2\le\Lambda\mathscr E_j\). Hence, we get \(\|q_j\|_2\le\Lambda\)
and \(\mathscr E_j\le\Lambda\).

Using the identity~\eqref{eq:B1} and \(|Q|=1\), we obtain
\[
\|\nabla u_j\|_2^2\ge\|\partial_j u_j\|_2^2\ge\left(\int_Q\partial_j u_j\dx\right)^2=1.
\]
Hence,
\[
\mathscr E_j\ge\lambda\|\nabla u_j\|_2^2\ge\lambda.
\]
This gives all the estimates in~\eqref{eq:B2}.

It remains to prove the maximum bound, which is adopted from \textsc{Stampacchia}'s truncation argument~\cite{SJL_1963-1964___3_1_0}.
\[
\zeta_+:=(u_1-1)^+\in H^1(Q),\qquad
\zeta_-:=u_1^-:=\max(-u_1,0)\in H^1(Q).
\]
A direct calculation gives
\[
\nabla\zeta_+(x)=\mathbf 1_{\{u_1>1\}}(x)\nabla u_1(x),\qquad
\nabla\zeta_-(x)=-\mathbf 1_{\{u_1<0\}}(x)\nabla u_1(x)
\quad\text{for a.e. }x\in Q.
\]
Since \(w\in V_1\) and \(u_1=x_1+w\), the trace of \(u_1\) equals \(0\) on \(\Gamma_L\) and \(1\) on \(\Gamma_R\). For almost every \(x_2\in(0,1)\), the absolutely continuous representative of \(u_1(\,\cdot\,,x_2)\) satisfies \(u_1(0,x_2)=0\) and \(u_1(1,x_2)=1\).
Composing this representative with the Lipschitz functions \(t\mapsto(t-1)^+\) and \(t\mapsto t^-\) shows that the corresponding slices of \(\zeta_+\) and \(\zeta_-\) vanish at both endpoints. Hence \(\zeta_+,\zeta_-\in V_1\).%

Taking \(\phi=\zeta_+\) in~\eqref{eq:W} and using the first gradient formula yields
\begin{align*}
0&=\int_{\{x\in Q:u_1(x)>1\}}b(x)\nabla u_1(x)\cdot\nabla u_1(x)\,dx\\
&=\int_Q b(x)\nabla\zeta_+(x)\cdot\nabla\zeta_+(x)\,dx\\
&\ge\lambda\|\nabla\zeta_+\|_{L^2(Q)}^2.
\end{align*}
Thus \(\nabla\zeta_+=0\) almost everywhere, and \Cref{lem:mixed_space_poincare} gives \(\zeta_+=0\) almost everywhere. Therefore \(u_1\le1\) almost everywhere. Similarly, testing~\eqref{eq:W} with \(\zeta_-\), multiplying by \(-1\), and repeating the preceding calculation gives \(\lambda\|\nabla\zeta_-\|_{L^2(Q)}^2\le0\). Another application of \Cref{lem:mixed_space_poincare} gives \(\zeta_-=0\) almost everywhere. This proves \(0\le u_1\le1\). 

Interchanging \(x_1,V_1\) with \(x_2,V_2\) in the same truncation argument gives \(0\le u_2\le1\).

\end{proof}

\subsection{Uniform invertibility of the mean-gradient matrix}
\label{sec:invertibility}
\begin{lemma}
\label{lem:G_invertible}
The matrix \(G\) defined in~\eqref{eq:FG} satisfies \(\sigma_{\min}(G)\ge c_*>0\) with \(\|G^{-1}\|_2\le c_*^{-1}\).
\end{lemma}

\begin{proof}
A direct calculation gives
\begin{equation}
G=\begin{pmatrix}1&\beta\\\gamma&1\end{pmatrix}\label{eq:G1}
\end{equation}
with
\[
\gamma=\int_0^1[u_1(r,1)-u_1(r,0)]\,dr,
\quad\text{and}\quad
\beta=\int_0^1[u_2(1,r)-u_2(0,r)]\,dr.
\]

We first prove \(|\gamma|\le1\) and \(|\beta|\le1\). By \Cref{lem:wellposed_energy}, \(0\le u_j\le1\) almost everywhere in \(Q\) for \(j=1,2\).
By \Cref{lem:rectangle_traces}, for almost every \(r\in(0,1)\), the slices \(u_1(r,\cdot)\) and \(u_2(\cdot,r)\) have absolutely continuous representatives whose endpoint values agree with the corresponding boundary traces. Since these representatives are continuous and take values in \([0,1]\) almost everywhere, their endpoint values also lie in \([0,1]\). Hence \(|u_1(r,1)-u_1(r,0)|\le1\) and \(|u_2(1,r)-u_2(0,r)|\le1\) for almost every \(r\); integration and \eqref{eq:G1} yield \(|\gamma|\le1\) and \(|\beta|\le1\).

We next obtain a uniform gap from equality. Let 
\[
t(r)=\operatorname{Tr}u_1(r,1),\quad b_0(r)=\operatorname{Tr}u_1(r,0),\quad\text{and\quad}
d=1-\gamma.
\]
Then
\[
d=\int_0^1[1-t(r)+b_0(r)]\,dr,
\]
so \(\int_0^1[1-t(r)]\,dr\le d\). We now prove the angular estimate used below. Fix \(0<\rho<1/2\), and parametrize the annular quarter-sector
\[
D_\rho:=\{(r\cos\theta,1-r\sin\theta):
\rho<r<1/2,\ 0<\theta<\pi/2\}
\]
by setting \(U(r,\theta):=u_1(r\cos\theta,1-r\sin\theta)\).
Because \(r\ge\rho>0\), this parametrization is a \(C^1\)
bi-Lipschitz diffeomorphism between
\((\rho,1/2)\times(0,\pi/2)\) and \(D_\rho\). 
This immediately implies \(U(r,\cdot)\in H^1(0,\pi/2)\) for almost every \(r\), with
\[
\partial_\theta U(r,\theta)
=-r\sin\theta\,\partial_1u_1
  -r\cos\theta\,\partial_2u_1
\quad\text{at }(r\cos\theta,1-r\sin\theta).
\]
In particular, \(|\partial_\theta U(r,\theta)|\le r|\nabla u_1(r\cos\theta,1-r\sin\theta)|\) almost everywhere.

Compatibility of Sobolev traces with bi-Lipschitz pullbacks yields \(U(r,0)=t(r)\) and \(U(r,\pi/2)=0\) for almost every \(r\in(\rho,1/2)\), where the second identity follows from the zero Dirichlet trace of \(u_1\) on \(\Gamma_L\). The fundamental theorem of calculus and the Cauchy--Schwarz inequality give
\[
t^2(r)
=\left|\int_0^{\pi/2}\partial_\theta U(r,\theta)\,d\theta\right|^2
\le\frac\pi2\int_0^{\pi/2}|\partial_\theta U(r,\theta)|^2\,d\theta.
\]
Using the polar Jacobian \(dx=r\,dr\,d\theta\)
and the preceding two inequalities, we obtain
\[
\begin{aligned}
\int_Q|\nabla u_1(x)|^2\,dx
&\ge \int_{D_\rho}|\nabla u_1(x)|^2\,dx\\
&=\int_\rho^{1/2}\int_0^{\pi/2}
 |\nabla u_1(r\cos\theta,1-r\sin\theta)|^2
 r\,d\theta\,dr\\
&\ge\int_\rho^{1/2}\frac1r
 \int_0^{\pi/2}|\partial_\theta U(r,\theta)|^2
 \,d\theta\,dr\\
&\ge\frac2\pi\int_\rho^{1/2}\frac{t^2(r)}{r}\,dr.
\end{aligned}
\]
Letting \(\rho\to 0\) and applying monotone convergence to the
nonnegative right-hand side proves
\begin{equation}
\int_Q|\nabla u_1(x)|^2\,dx
\ge\frac2\pi\int_0^{1/2}\frac{t^2(r)}{r}\,dr.
\label{eq:G2}
\end{equation}

If \(d=0\), then \(t=1\) almost everywhere and the right-hand side of~\eqref{eq:G2} diverges, a contradiction. If \(d\ge1/4\), then trivially \(d\ge c_*\). It remains to consider \(0<d<1/4\). Since \(0\le t\le1\) and \(t^2\ge1-2(1-t)\),
\[
\begin{aligned}
\int_{2d}^{1/2}\frac{t^2(r)}{r}\,dr
&\ge\log\frac1{4d}
-2\int_{2d}^{1/2}\frac{1-t(r)}r\,dr\\
&\ge\log\frac1{4d}-\frac1d\int_0^1(1-t(r))\,dr\\
&\ge\log\frac1{4d}-1.
\end{aligned}
\]
Combining this estimate with~\eqref{eq:G2} and~\eqref{eq:B2}, we obtain \(\log(1/(4d))-1\le\pi\Lambda/(2\lambda)\), and hence \(1-\gamma=d\ge c_*\).

To check the opposite side without suppressing signs, let \(d'=1+\gamma=\int_0^1[t(r)+1-b_0(r)]\,dr\). Then \(\int_0^1[1-b_0(r)]\,dr\le d'\). Repeating the same estimate in the lower-left sector, with \(b_0\) in place of \(t\), gives \(1+\gamma\ge c_*\). Therefore \(|\gamma|\le1-c_*\).

For the second loading, we employ the coordinate-swap matrix 
\[
P=\begin{pmatrix}0&1\\1&0\end{pmatrix}.
\]
The function \(\widehat u(x)=u_2(Px)\) satisfies the first-loading problem with coefficient \(P b(Px)P^{\top}\), whose ellipticity parameters are unchanged. Its off-diagonal mean is exactly \(\beta\). Thus \(|\beta|\le1-c_*\).

Finally, write \(G=I+K\). By \eqref{eq:G1}, \(K^TK=\operatorname{diag}(\gamma^2,\beta^2)\), and hence \(\|K\|_2\le1-c_*\). Therefore \(|G\xi|\ge|\xi|-|K\xi|\ge c_*|\xi|\) for every vector \(\xi\), so \(\sigma_{\min}(G)\ge c_*\). Consequently, \(G\) is invertible and \(\|G^{-1}\|_2\le c_*^{-1}\), which completes the proof.
\end{proof}
\subsection{Direct construction of the stream function and the dual mixed problem}\label{sec:stream}

In this part, we prove two cross-gradient orthogonality identities, which are key to the final error estimate. By~\cite[\S 3, Theorem 1.1]{Necas2012}, an integration by parts, we have
\begin{lemma}
\label{lem:cross_gradient}
If \(w\in V_1\) and \(\psi\in V_2\), then
\[
\int_QJ\nabla w(x)\cdot\nabla\psi(x)\,dx
=\int_QJ\nabla\psi(x)\cdot\nabla w(x)\,dx=0.
\]
\end{lemma}

The following lemma constructs a stream function and identifies the associated dual mixed problem. The construction agrees with the standard stream-function relation for planar divergence-form equations~\cite[\S1.3, equations~(1.15)--(1.16)]{AlessandriniNesi2009}, while additionally proving the complete exchange of the mixed boundary conditions in the present setting.
\begin{lemma}\label{lem:stream}
Let \(u\) be the first-loading weak solution for the coefficient \(b\), let \(q=b\nabla u\), and let \(\mathscr E\) denote its energy. There exists \(s\in H^1(Q)\) such that \(\nabla s=Jq\), \(s-\mathscr E x_2\in V_2\), and
\begin{equation*}
\int_Qb^\#(x)\nabla s(x)\cdot\nabla\psi(x)\dx=0,
\quad\text{for every }\psi\in V_2,
\end{equation*}
with \(0\le s\le\mathscr E\). 
\end{lemma}

\begin{proof}
Define
\begin{equation}
s(x_1,x_2){:}=\int_0^{x_2}e_1\cdot q(x_1,t)\,\mathrm{d}t.
\label{eq:S0}
\end{equation}
Fubini's theorem and the Cauchy--Schwarz inequality give \(s\in L^2(Q)\) and \(\partial_2s=e_1\cdot q\). To identify $\pa_1 s$, we fix an arbitrary \(\phi\in C_c^\infty(Q)\) and set
\[
\Psi(x_1,x_2):=\int_{x_2}^1\phi(x_1,y)\,\mathrm{d}y.
\]
The zero extension of \(\phi\) is smooth, and hence
\(\Psi\in C^\infty(\overline Q)\subset H^1(Q)\). Because \(\phi\) is
compactly supported in \(Q\), it vanishes for \(x_1\) near \(0\) and
\(1\). Hence \(\Psi\) has zero trace on \(\Gamma_L\cup\Gamma_R\), so
\(\Psi\in V_1\). Moreover,
\[
\partial_2\Psi(x_1,x_2)=-\phi(x_1,x_2),
\qquad
\partial_1\Psi(x_1,x_2)
=\int_{x_2}^1\partial_1\phi(x_1,y)\,\mathrm{d}y.
\]

We first justify the Fubini interchange. Since \(e_1\cdot q\in L^2(Q)\)
and \(\abs{Q}=1\), we obtain
\[
\int_0^1\int_0^1\int_0^{x_2}
 |e_1\cdot q(x_1,t)|\,|\partial_1\phi(x_1,x_2)|
 \,\mathrm{d}t\dx_2\dx_1\le
\|\partial_1\phi\|_{L^\infty(Q)}
\|e_1\cdot q\|_{L^1(Q)}<\infty.
\]
Therefore, we may apply Fubini's theorem on the triangular region
\(0<t<x_2<1\), and use~\eqref{eq:S0} to get
\[
\begin{aligned}
\int_Qs(x)\,\partial_1\phi(x)\dx
&=\int_0^1\int_0^1
 \left(\int_0^{x_2}e_1\cdot q(x_1,t)\,\mathrm{d}t\right)
 \partial_1\phi(x_1,x_2)\dx_2\dx_1\\
&=\int_0^1\int_0^1
 (e_1\cdot q(x_1,t))
 \left(\int_t^1\partial_1\phi(x_1,x_2)\dx_2\right)
 \,\mathrm{d}t\dx_1\\
&=\int_Q(e_1\cdot q(x))\partial_1\Psi(x)\dx.
\end{aligned}
\]
Since \(\Psi\in V_1\), it is admissible. Using \(\partial_2\Psi=-\phi\), we obtain
\[
\begin{aligned}
0
&=\int_Qq(x)\cdot\nabla\Psi(x)\,dx\\
&=\int_Q(e_1\cdot q(x))\partial_1\Psi(x)\,dx
  -\int_Q(e_2\cdot q(x))\phi(x)\,dx.
\end{aligned}
\]
Combining the preceding two identities, we obtain, for every
\(\phi\in C_c^\infty(Q)\),
\[
\int_Qs(x)\partial_1\phi(x)\,dx
=\int_Q(e_2\cdot q(x))\phi(x)\,dx.
\]
This gives
\(\partial_1s=-e_2\cdot q\). Together with
\(\partial_2s=e_1\cdot q\), this yields \(\nabla s=Jq\), and hence
\(s\in H^1(Q)\).

By Lemma~\ref{lem:rectangle_traces}, the traces of \(s\) agree with
the endpoint values of its almost-everywhere vertical slices. Formula~\eqref{eq:S0}
gives the bottom endpoint \(0\) and the top endpoint
\(c(x_1)=\int_0^1 e_1\cdot q(x_1,t)\,dt\). By the Cauchy--Schwarz inequality,
\(\|c\|_{L^2(0,1)}\le\|q\|_{L^2(Q)}\), so \(c\in L^2(0,1)\).
The trace claims follow from the two estimates
\[
\begin{aligned}
\|s(\cdot,h)\|_{L^2(0,1)}^2
&\le h\int_0^h\|e_1\cdot q(\cdot,t)\|_{L^2(0,1)}^2\,\mathrm{d}t\to0,\\
\|s(\cdot,1-h)-c\|_{L^2(0,1)}^2
&\le h\int_{1-h}^1\|e_1\cdot q(\cdot,t)\|_{L^2(0,1)}^2\,\mathrm{d}t\to0.
\end{aligned}
\]
We now determine the top trace \(c\). Fix an arbitrary
\(\zeta\in C_c^\infty(0,1)\) and define
\(Z(x_1,x_2):=\zeta(x_1)\) on \(Q\). Then
\(Z\in C^\infty(\overline Q)\subset H^1(Q)\). Because \(\zeta\) is
compactly supported in \((0,1)\),
\(\operatorname{Tr}_{\Gamma_L}Z=\zeta(0)=0\) and
\(\operatorname{Tr}_{\Gamma_R}Z=\zeta(1)=0\).
Thus \(Z\in V_1\), and \(\nabla Z=(\zeta'(x_1),0)^T\). Since
\((e_1\cdot q)\zeta'\in L^1(Q)\), using the first-loading weak formulation
\eqref{eq:W} and Fubini's theorem, we get
\begin{equation*}
\begin{aligned}
0
&=\int_Qq(x)\cdot\nabla Z(x)\,dx\\
&=\int_0^1\int_0^1
   \bigl(e_1\cdot q(x_1,x_2)\bigr)\zeta'(x_1)\,dx_2\,dx_1\\
&=\int_0^1c(x_1)\zeta'(x_1)\,dx_1.
\end{aligned}
\end{equation*}
Since \(\zeta\) is arbitrary, \(c'=0\) in \(\mathscr D'(0,1)\).
This immediately implies \(c=C\) in \(\mathscr D'(0,1)\) for some constant \(C\).
Since \(c\in L^2(0,1)\), we have \(c(x_1)=C\) for almost every
\(x_1\in(0,1)\). Integrating over \((0,1)\) and testing the mixed weak equation with \(u-x_1\in V_1\) gives \(C=\int_Q e_1\cdot q\,dx=\mathscr E\).
Consequently, \(\operatorname{Tr}_{\Gamma_B}s=0\) and
\(\operatorname{Tr}_{\Gamma_T}s=\mathscr E\).
Hence \(s-\mathscr E x_2\in V_2\).

\Needspace{8\baselineskip}
Lemma~\ref{lem:matrix_tools} and \(J^2=-I\) give \(b^\#\nabla s=J\nabla u\) almost everywhere. Write \(u=x_1+w\) with \(w\in V_1\). For every \(\psi\in V_2\), we use Lemma~\ref{lem:cross_gradient} to get
\[
\begin{aligned}
\int_Qb^\#(x)\nabla s(x)\cdot\nabla\psi(x)\dx
&=\int_QJe_1\cdot\nabla\psi(x)\dx
  +\int_QJ\nabla w(x)\cdot\nabla\psi(x)\dx\\
&=\int_Q\partial_2\psi(x)\dx=0.
\end{aligned}
\]
This proves the complete dual weak formulation.
Finally, \(b^\#\) is uniformly elliptic and \(\mathscr E>0\). 

Testing the dual weak formulation with \((s-\mathscr E)^+\) and \(s^-\), respectively, and repeating the truncation argument from \Cref{lem:wellposed_energy} yields \(0\le s\le\mathscr E\), which completes the proof.
\end{proof}
\subsection{Weak boundary functionals and boundary measures}
\label{sec:boundary-functionals}
In what follows, we interpret the boundary fluxes and tangential derivatives as weak boundary functionals and prove that they are finite signed Radon measures.
\begin{definition}
\label{def:boundary_functionals}
For the \(j\)-th loading, let \(q_j=b\nabla u_j\) and define
\[
\Gamma_D^1{:}=\Gamma_L\cup\Gamma_R,\qquad
\Gamma_N^1{:}=\Gamma_B\cup\Gamma_T,\qquad
\Gamma_D^2{:}=\Gamma_N^1,\qquad
\Gamma_N^2{:}=\Gamma_D^1.
\]

If \(r\in L^2(Q;\mathbb R^2)\) and \(\operatorname{div}r=0\) in the sense of distributions, define its normal trace by
\begin{equation}
\langle\gamma_n r,\operatorname{Tr}\Phi\rangle_{\partial Q}
{:}=\int_Q r(x)\cdot\nabla\Phi(x)\dx,
\qquad \text{for all\quad}\Phi\in H^1(Q).
\label{eq:NT}
\end{equation}

The weak formulation gives \(\gamma_nq_j=0\) on \(\Gamma_N^j\). Moreover, \(\operatorname{div}(J\nabla u_j)=0\), and the oriented boundary tangent \(\tau=Jn\), where \(n\) is the outward unit normal on the edges of \(Q\), satisfies the identity
\begin{equation}
D_\tau(\operatorname{Tr}u_j)=-\gamma_n(J\nabla u_j)
\label{eq:TN}
\end{equation}
in the sense of distributions. Since \(u_j\) is constant on each connected component of \(\Gamma_D^j\), the right-hand side of~\eqref{eq:TN} vanishes on \(\Gamma_D^j\). We write
\[
\nu_j=\gamma_nq_j\big|_{\Gamma_D^j},
\qquad
\mu_j=D_\tau(\operatorname{Tr}u_j)\big|_{\Gamma_N^j}.
\]
\end{definition}

\begin{figure}[htbp]
\centering
\begin{tikzpicture}[
    x=1cm,y=1cm,
    >={Latex[length=2.4mm,width=1.7mm]},
    font=\small,
    dirichlet/.style={draw=dnBlue,line width=1.8pt,line cap=round},
    neumann/.style={draw=dnOrange,line width=1.8pt,
                    dash pattern=on 5pt off 3pt,line cap=round},
    vector/.style={->,draw=black!85,line width=0.95pt},
    edge-label/.style={font=\footnotesize,inner sep=1pt}
  ]
  \definecolor{dnBlue}{RGB}{0,114,178}
  \definecolor{dnOrange}{RGB}{213,94,0}

  \draw[dirichlet] (0.18,4.43)--(0.73,4.43);
  \node[anchor=west,font=\footnotesize] at (0.88,4.43)
    {Dirichlet: $n$, $\nu_j$};
  \draw[neumann] (5.00,4.43)--(5.55,4.43);
  \node[anchor=west,font=\footnotesize] at (5.70,4.43)
    {Neumann: $\tau=Jn$, $\mu_j$};

  \node[font=\bfseries] at (1.50,3.80) {Loading $j=1$};
  \fill[black!2] (0,0) rectangle (3,3);
  \draw[black!15,line width=0.35pt] (0,0) rectangle (3,3);
  \draw[dirichlet] (0,0)--(0,3);
  \draw[dirichlet] (3,0)--(3,3);
  \draw[neumann] (0,0)--(3,0);
  \draw[neumann] (0,3)--(3,3);
  \node[font=\normalsize] at (1.50,1.50) {$Q$};
  \node[edge-label,anchor=east,text=dnBlue] at (-0.18,1.50) {$\Gamma_L$};
  \node[edge-label,anchor=west,text=dnBlue] at (3.18,1.50) {$\Gamma_R$};
  \node[edge-label,anchor=south,text=dnOrange] at (1.50,3.15) {$\Gamma_T$};
  \node[edge-label,anchor=north,text=dnOrange] at (1.50,-0.15) {$\Gamma_B$};

  \draw[vector] (0,2.08)--(-0.49,2.08);
  \draw[vector] (3,2.08)--(3.49,2.08);
  \draw[vector] (2.09,2.84)--(1.13,2.84);
  \draw[vector] (1.04,0.16)--(2.00,0.16);

  \begin{scope}[shift={(6.10,0)}]
    \node[font=\bfseries] at (1.50,3.80) {Loading $j=2$};
    \fill[black!2] (0,0) rectangle (3,3);
    \draw[black!15,line width=0.35pt] (0,0) rectangle (3,3);
    \draw[neumann] (0,0)--(0,3);
    \draw[neumann] (3,0)--(3,3);
    \draw[dirichlet] (0,0)--(3,0);
    \draw[dirichlet] (0,3)--(3,3);
    \node[font=\normalsize] at (1.50,1.50) {$Q$};

    \draw[vector] (1.50,3)--(1.50,3.50);
    \draw[vector] (1.50,0)--(1.50,-0.50);
    \draw[vector] (0.16,2.35)--(0.16,1.39);
    \draw[vector] (2.84,1.39)--(2.84,2.35);

  \end{scope}
\end{tikzpicture}
\caption{Boundary partitions for the two coordinate loadings:
Dirichlet (solid blue) and homogeneous Neumann (dashed vermilion) edges.
Arrows indicate the outward normal \(n\) on \(\Gamma_D^j\) and the oriented
tangent \(\tau=Jn\) on \(\Gamma_N^j\); \(\nu_j\) and \(\mu_j\) are the
corresponding weak boundary functionals.}
\label{fig:mixed-boundary-measures}
\end{figure}
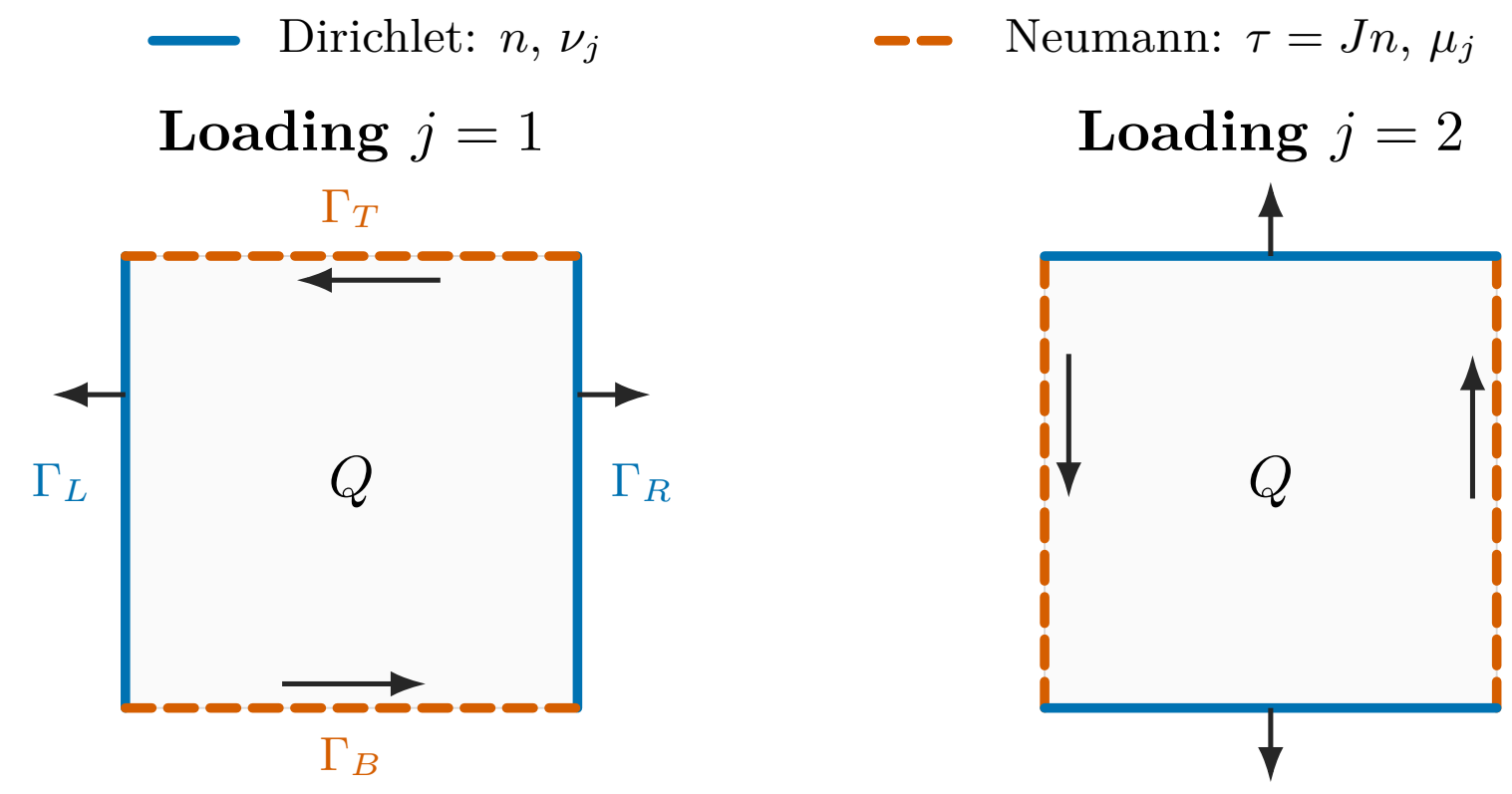

The boundary assignments and orientations are illustrated in
\Cref{fig:mixed-boundary-measures}.
\FloatBarrier

In what follows, we discuss the total variations of the boundary measures.
\begin{lemma}
\label{lem:boundary_measures}
For the first loading, write \(u=u_1\), \(q=q_1\), and
\(\mathscr E=\int_Qq(x)\cdot\nabla u(x)\,dx\). The boundary functionals \(\nu=\nu_1\) and \(\mu=\mu_1\) are finite signed Radon measures and satisfy
\begin{equation}
\begin{aligned}
&\nu\!\lfloor_{\Gamma_L}\le0,\qquad
 \nu\!\lfloor_{\Gamma_R}\ge0,\qquad
 |\nu|(\Gamma_L\cup\Gamma_R)=2\mathscr E\le2\Lambda,\\
&\mu\!\lfloor_{\Gamma_B}\ge0,\qquad
 \mu\!\lfloor_{\Gamma_T}\le0,\qquad
 |\mu|(\Gamma_B\cup\Gamma_T)=2.
\end{aligned}
\label{eq:BV}
\end{equation}
After interchanging the two pairs of edges, the second loading satisfies the same total-variation bounds, with \(\mathscr E\) in the first line replaced by the energy of the second loading.
\end{lemma}

\begin{proof}
We first prove a sign principle for weak conormal functionals and then apply it to the four relevant edges.

\medskip\noindent\textbf{Step 1. Sign and measure property of a conormal functional.}
Let \(0<\alpha\le\beta\), \(c\in\mathcal M(\alpha,\beta;Q)\), and let \(E,E'\) be the two closed Dirichlet edges of a mixed problem. Set \(V_D:=\{\psi\in H^1(Q):\operatorname{Tr}\psi=0\text{ on }E\cup E'\}\). Suppose that \(v\in H^1(Q)\) satisfies \(\int_Qc\nabla v\cdot\nabla\psi\dx=0\) for every \(\psi\in V_D\), that \(v\ge0\) in \(Q\), and that \(\operatorname{Tr}v=0\) on \(E\). For \(\varphi\in C^\infty(E)\), choose a smooth lifting \(\Phi\) which equals \(\varphi\) on \(E\) and vanishes on \(E'\), and define
\[
L_E(\varphi){:}=\int_Qc\nabla v\cdot\nabla\Phi\dx.
\]
The difference of two such liftings belongs to \(V_D\), so this value is independent of the lifting. If \(\varphi\ge0\), the lifting may be chosen nonnegative. Define \(T_\varepsilon(t):=\min\{1,\max\{0,t/\varepsilon\}\}\).
The function \(\Phi T_\varepsilon(v)\) belongs to \(V_D\): its trace vanishes on \(E\) because \(T_\varepsilon(v)=0\) there, and on \(E'\) because \(\Phi=0\) there.
Testing with it gives \(0=I_\varepsilon+J_\varepsilon\), where
\begin{align*}
I_\varepsilon&{:}=\int_QT_\varepsilon(v)c\nabla v\cdot\nabla\Phi\dx,\\
J_\varepsilon&{:}=\frac{1}{\varepsilon}\int_{\{0<v<\varepsilon\}}\Phi\,\nabla v\cdot c\nabla v\dx.
\end{align*}
Since \(\Phi\ge0\) and \(\nabla v\cdot c\nabla v\ge\alpha|\nabla v|^2\), we have \(J_\varepsilon\ge0\) and hence \(I_\varepsilon\le0\). Moreover, \(\nabla v=0\) almost everywhere on \(\{v=0\}\)~\cite[Theorem~4.4]{EvansGariepy2015}, so \(T_\varepsilon(v)c\nabla v\to c\nabla v\) almost everywhere in \(Q\). The bound \(|c\nabla v|\le\beta|\nabla v|\) and dominated convergence in \(L^2\) yield \(I_\varepsilon\to L_E(\varphi)\). Consequently, \(L_E(\varphi)\le0\) for every \(\varphi\ge0\). If instead \(M-v\ge0\) and \(\operatorname{Tr}v=M\) on \(E\), applying the same argument to \(M-v\) gives \(L_E(\varphi)\ge0\).

It remains to justify that this sign-definite functional is a measure. If \(L_E\le0\), then \(-L_E\) is positive. Monotonicity and \(-\|\varphi\|_\infty\le\varphi\le\|\varphi\|_\infty\) imply \(|L_E(\varphi)|\le-L_E(1)\|\varphi\|_\infty\). Since \(C^\infty(E)\) is uniformly dense in \(C(E)\), \(-L_E\) extends uniquely to a positive continuous functional on \(C(E)\). The Riesz theorem~\cite[Theorem~1.39]{EvansGariepy2015} represents it by a finite positive Radon measure, so \(L_E\) is a finite signed Radon measure. The case \(L_E\ge0\) is identical.

The representing measure has no atoms at the endpoints of \(E\). Indeed, fix an endpoint \(p\) and choose a smooth radial cutoff \(0\le\Psi_\rho\le1\), equal to one on \(B_{\rho^2}(p)\), supported in \(B_\rho(p)\), and with a smoothed logarithmic transition between these circles. Then \(\|\nabla\Psi_\rho\|_{L^2(Q)}^2\le C/|\log\rho|\). For small \(\rho\), \(\Psi_\rho\) vanishes on \(E'\) and is an admissible lifting. Since \(L_E\) has a fixed sign,
\[
|L_E|(\{p\})\le |L_E(\Psi_\rho|_E)|
\le\|c\nabla v\|_{L^2(Q)}\|\nabla\Psi_\rho\|_{L^2(Q)}\longrightarrow0.
\]
Thus the measure has the same total mass on \(E\) and on its open interior.

\medskip\noindent\textbf{Step 2. The flux measures.}
For \(\varphi\in C^\infty([0,1])\), define
\[
\begin{aligned}
\nu_R(\varphi)&{:}=\int_Qq\cdot\nabla\bigl(x_1\varphi(x_2)\bigr)\dx,\\
\nu_L(\varphi)&{:}=\int_Qq\cdot\nabla\bigl((1-x_1)\varphi(x_2)\bigr)\dx.
\end{aligned}
\]
For smooth test functions compactly supported on either open vertical edge, the corresponding affine lifting has zero trace on the other edges. Thus the restrictions of these measures to the open edges agree with \(\gamma_nq\) in \Cref{def:boundary_functionals}. By \Cref{lem:wellposed_energy}, \(0\le u\le1\). Step~1 applied to \(1-u\) on \(\overline{\Gamma_R}\), whose flux is \(-q\), gives \(\nu_R\ge0\); applied to \(u\) on \(\overline{\Gamma_L}\), it gives \(\nu_L\le0\). Taking \(\varphi\equiv1\) and testing~\eqref{eq:W} with \(u-x_1\in V_1\) gives \(\nu_R(1)=\mathscr E\) and \(\nu_L(1)=-\mathscr E\).
Since Step~1 excludes endpoint atoms, \(|\nu|(\Gamma_L\cup\Gamma_R)=2\mathscr E\le2\Lambda\).

\medskip\noindent\textbf{Step 3. The tangential-derivative measures.}
Let \(s\) be the stream function from Lemma~\ref{lem:stream} and set \(r=b^\#\nabla s=J\nabla u\).
The dual weak formulation has \(\Gamma_B\cup\Gamma_T\) as its Dirichlet boundary, and
\(0\le s\le\mathscr E\), with traces \(s=0\) on \(\Gamma_B\) and \(s=\mathscr E\) on \(\Gamma_T\). Define
\[
\begin{aligned}
\kappa_T(\varphi)&{:}=\int_Qr\cdot\nabla\bigl(x_2\varphi(x_1)\bigr)\dx,\\
\kappa_B(\varphi)&{:}=\int_Qr\cdot\nabla\bigl((1-x_2)\varphi(x_1)\bigr)\dx.
\end{aligned}
\]
Applying Step~1 to \(\mathscr E-s\) on \(\overline{\Gamma_T}\) and to \(s\) on \(\overline{\Gamma_B}\) yields \(\kappa_T\ge0\) and \(\kappa_B\le0\). For \(\varphi\equiv1\), the identity \(r=J\nabla u\) and the first identity in~\eqref{eq:B1} give \(\kappa_T(1)=1\) and \(\kappa_B(1)=-1\). By~\eqref{eq:TN},
\(\mu\!\lfloor_{\Gamma_T}=-\kappa_T\) and \(\mu\!\lfloor_{\Gamma_B}=-\kappa_B\). The preceding mass identities and the absence of endpoint atoms yield \(\mu\!\lfloor_{\Gamma_B}\ge0\), \(\mu\!\lfloor_{\Gamma_T}\le0\), and \(|\mu|(\Gamma_B\cup\Gamma_T)=2\).
This proves~\eqref{eq:BV} for the first loading. 

Since the orthogonal coordinate swap preserves the coefficient parameters and all preceding estimates, interchanging \(x_1,x_2\) and the two pairs of edges proves the second-loading statement, which completes the proof.
\end{proof}
\section{Error estimate and sharpness for mixed RVE reconstruction}
\label{sec:completion}
We henceforth take \(b=a_\eta\) in the preceding results and use the periodic correctors in~\eqref{eq:H}.
Let \(\chi_i^*\in H^1_{\rm per}(Q)\) be the zero-mean periodic corrector for \(a^{\top}\), so that \(\int_Qa^{\top}(e_i+\nabla\chi_i^*)\cdot\nabla\phi\,dy=0\) for all \(\phi\in H^1_{\rm per}(Q)\). The corresponding effective matrix is \(\mathcal A^{\top}\), by the adjoint corrector identity~\cite{Shen2018}.

The following characterization of the periodic flux potentials and uniform bounds is also well-known; we keep a short proof for the sake of self-evidence.
\begin{lemma}\label{lem:flux_potentials}
For some \(p>2\) depending only on \(\lambda\) and \(\Lambda\), there exist zero-mean periodic functions \(\theta_i^*\) such that \(\chi_i^*,\theta_i^*\in W^{1,p}_{\rm per}(Q)\) and
\begin{equation}
\begin{aligned}
h_i^*&=a^{\top}(e_i+\nabla\chi_i^*)-\mathcal A^{\top}e_i
=-J\nabla\theta_i^*,\\
\max_{i=1,2}&\big(\|\chi_i^*\|_{L^\infty(Q)}+
\|\theta_i^*\|_{L^\infty(Q)}\big)
\le C_{\rm cor}(\lambda,\Lambda).
\end{aligned}
\label{eq:C2}
\end{equation}
Their \(W^{1,p}(Q)\) norms are also bounded in terms of \(\lambda\) and \(\Lambda\).
\end{lemma}

\begin{proof}
We exploit the flux-corrector construction for the adjoint coefficient \(a^{\top}\) in~\cite{Shen2018}. In the present scalar case, it yields periodic functions \(\phi_{k\ell i}^*\in H^1_{\rm per}(Q)\), antisymmetric in \(k,\ell\), such that
\[
(h_i^*)_\ell=\partial_k\phi_{k\ell i}^*,
\qquad
\phi_{k\ell i}^*=-\phi_{\ell k i}^*.
\]
In two dimensions, set \(\theta_i^*:=-\phi_{12i}^*\) and subtract its mean. Antisymmetry then gives the flux-potential identity in~\eqref{eq:C2}.

By the periodic Meyers estimate~\cite{Meyers1963,Shen2018}, there is \(p>2\), depending only on \(\lambda\) and \(\Lambda\), for which \(\|\chi_i^*\|_{W^{1,p}(Q)}\le C\). Since \(a\) and \(\mathcal A\) are bounded in terms of \(\lambda\) and \(\Lambda\), \(h_i^*\in L^p(Q)\) with a uniform bound. The first identity in~\eqref{eq:C2} gives \(\nabla\theta_i^*=Jh_i^*\); the zero-mean Poincar\'e inequality therefore bounds \(\theta_i^*\) in \(W^{1,p}(Q)\). The two-dimensional Morrey embedding gives the stated \(L^\infty\) bounds.
\end{proof}
The following identity pairs the boundary measures with continuous periodic correctors.
\begin{lemma}
\label{lem:boundary_identity}
Assume only the bounded, measurable, and periodic coefficient hypotheses of \Cref{sec:problem}. For the \(j\)-th mixed loading, let \(q_j=a_\eta\nabla u_j\), and let \(\nu_j,\mu_j\) be the boundary measures from \Cref{def:boundary_functionals} and \Cref{lem:boundary_measures}. For a periodic function \(f\), write \(f^\eta(x)=f(\vartheta+x/\eta)\). Then, for every \(i,j\in\{1,2\}\),
\begin{equation}
e_i\cdot\int_Q\bigl(q_j-\mathcal A\nabla u_j\bigr)\,dx
=-\eta\left[
\int_{\Gamma_D^j}\chi_i^{*,\eta}\,d\nu_j
+\int_{\Gamma_N^j}\theta_i^{*,\eta}\,d\mu_j
\right].
\label{eq:R}
\end{equation}
\end{lemma}

The right-hand side of~\eqref{eq:R} pairs continuous functions with finite Radon measures.
By Lemma~\ref{lem:flux_potentials}, both \(\chi_i^{*,\eta}\) and \(\theta_i^{*,\eta}\) are continuous on \(\overline Q\); hence the identity is well defined for measurable coefficients.

\begin{proof}
We recall two Green formulas:
If \(w\in W^{1,p}_{\rm per}(Q)\) with \(p>2\), then
\begin{equation}
\begin{aligned}
\int_Qq_j\cdot\nabla w^\eta\,dx
&=\int_{\Gamma_D^j}w^\eta\,d\nu_j,\\
\int_QJ\nabla u_j\cdot\nabla w^\eta\,dx
&=-\int_{\Gamma_N^j}w^\eta\,d\mu_j.
\end{aligned}
\label{eq:weak-green-multipliers}
\end{equation}
Choose smooth periodic \(w_m\) such that
\(w_m\to w\) in \(W^{1,p}_{\rm per}(Q)\). By the continuous Morrey embedding~\cite[Theorem~4.12, Part~II]{AdamsFournier2003}, the continuous periodic representatives converge uniformly. Since \(\eta\) is fixed, \(w_m^\eta\to w^\eta\) in \(W^{1,p}(Q)\) and uniformly on \(\overline Q\).

Fix \(j\in\{1,2\}\), set \(k=3-j\), and define the affine liftings
\[
\begin{aligned}
L_{D,m}^j(x)
&:=(1-x_j)w_m^\eta(x-x_je_j)
  +x_jw_m^\eta\bigl(x+(1-x_j)e_j\bigr),\\
L_{N,m}^j(x)
&:=(1-x_k)w_m^\eta(x-x_ke_k)
  +x_kw_m^\eta\bigl(x+(1-x_k)e_k\bigr).
\end{aligned}
\]
The first lifting agrees with \(w_m^\eta\) on \(\Gamma_D^j\), whereas the second agrees with it on \(\Gamma_N^j\). Therefore \(w_m^\eta-L_{D,m}^j\in V_j\) and \(w_m^\eta-L_{N,m}^j\in V_k\).
The weak equation for \(u_j\) and the affine-lifting representation of \(\nu_j\) in the proof of \Cref{lem:boundary_measures} give
\[
\int_Qq_j\cdot\nabla w_m^\eta\,dx
=\int_Qq_j\cdot\nabla L_{D,m}^j\,dx
=\int_{\Gamma_D^j}w_m^\eta\,d\nu_j.
\]

For the second identity, first note that \(u_j-x_j\in V_j\). By Lemma~\ref{lem:cross_gradient}, \(\int_QJ\nabla(u_j-x_j)\cdot\nabla\phi\,dx=0\) for every \(\phi\in V_k\). Also \(Je_j=(-1)^{j+1}e_k\), and the trace of \(\phi\) vanishes on the two edges normal to \(e_k\). Therefore \(\int_QJ\nabla u_j\cdot\nabla\phi\,dx=0\) for every \(\phi\in V_k\).
Taking \(\phi=w_m^\eta-L_{N,m}^j\) and using
\(\gamma_n(J\nabla u_j)=-\mu_j\) on \(\Gamma_N^j\), we obtain
\[
\int_QJ\nabla u_j\cdot\nabla w_m^\eta\,dx
=\int_QJ\nabla u_j\cdot\nabla L_{N,m}^j\,dx
=-\int_{\Gamma_N^j}w_m^\eta\,d\mu_j.
\]
Letting \(m\to\infty\) proves~\eqref{eq:weak-green-multipliers}. The volume terms converge because
\(q_j,J\nabla u_j\in L^2(Q)\) and
\(\nabla w_m^\eta\to\nabla w^\eta\) in
\(L^p(Q)\subset L^2(Q)\); the boundary terms converge by uniform convergence and the finiteness of the measures.

By Lemma~\ref{lem:flux_potentials}, for some \(p>2\),
\(\chi_i^*,\theta_i^*\in W^{1,p}_{\rm per}(Q)\), and the two representations of \(h_i^*\) are given by~\eqref{eq:C2}.
The chain rule under rescaling gives
\((\nabla_y\chi_i^*)^\eta=\eta\nabla_x\chi_i^{*,\eta}\) and
\(h_i^{*,\eta}=-\eta J\nabla_x\theta_i^{*,\eta}\).
Taking the scalar product with \(\nabla u_j\) in the two representations of \(h_i^{*,\eta}\) yields
\begin{align*}
h_i^{*,\eta}\cdot\nabla u_j
&=e_i\cdot(q_j-\mathcal A\nabla u_j)
+(\nabla_y\chi_i^*)^\eta\cdot q_j,\\
h_i^{*,\eta}\cdot\nabla u_j
&=\eta\nabla_x\theta_i^{*,\eta}\cdot J\nabla u_j.
\end{align*}
Both equalities hold in \(L^r(Q)\), where
\(r=2p/(p+2)>1\). Indeed, the corrector gradients and \(h_i^*\) lie in \(L^p\), while \(q_j\) and \(\nabla u_j\) lie in \(L^2\). Subtracting the two identities and integrating gives
\[
e_i\cdot\int_Q(q_j-\mathcal A\nabla u_j)\,dx
=\eta\left[\int_Q\nabla\theta_i^{*,\eta}\cdot J\nabla u_j\,dx
-\int_Qq_j\cdot\nabla\chi_i^{*,\eta}\,dx\right].
\]
Applying~\eqref{eq:weak-green-multipliers} first to \(w=\theta_i^*\) and then to \(w=\chi_i^*\) transforms the right-hand side into the boundary expression in~\eqref{eq:R}, which completes the proof.
\end{proof}

We are ready to prove the main result.
\begin{proof}[Proof of \Cref{thm:main}]
Equation~\eqref{eq:R}, Lemma~\ref{lem:boundary_measures}, and the estimate for finite measures give
\begin{equation}
\begin{aligned}
\left|e_i\cdot\int_Q
\big(q_j(x)-\mathcal A\nabla u_j(x)\big)\,dx\right|
&\le\eta\left(
\|\chi_i^*\|_\infty|\nu_j|(\Gamma_D^j)
+\|\theta_i^*\|_\infty|\mu_j|(\Gamma_N^j)\right)\\
&\le2\eta\left(
\Lambda\|\chi_i^*\|_\infty+\|\theta_i^*\|_\infty\right).
\end{aligned}
\label{eq:ENTRY}
\end{equation}

We use equations~\eqref{eq:C2} and~\eqref{eq:ENTRY} to control all four matrix entries, and therefore
\begin{equation}
\|F_\eta-\mathcal A G_\eta\|_2
\le\|F_\eta-\mathcal A G_\eta\|_F
\le4\eta(\Lambda+1)C_{\rm cor}.
\label{eq:RES}
\end{equation}
By Lemma~\ref{lem:G_invertible}, \(G_\eta\) is invertible. We now use
\begin{equation}
\mathcal A_\eta-\mathcal A
=(F_\eta-\mathcal A G_\eta)G_\eta^{-1}.
\label{eq:REC}
\end{equation}
A combination of~\eqref{eq:RES}--\eqref{eq:REC} with \(\|G_\eta^{-1}\|_2\le c_*^{-1}\) gives
\[
\|\mathcal A_\eta-\mathcal A\|_2
\le4c_*^{-1}(\Lambda+1)C_{\rm cor}\,\eta.
\]
This proves the error estimate in \Cref{thm:main}. The constant depends only on \(\lambda\) and \(\Lambda\) and is uniform in the phase \(\vartheta\). 
\end{proof}

A simple change of variables leads to \Cref{cor:physical_rate}.
\begin{proof}
Set \(z=(x-x_*)/\delta\), \(\eta=\epsilon/\delta\), \(\vartheta=x_*/\epsilon\), and \(u_j(z)=[v_j^\epsilon(x_*+\delta z)-x_{*,j}]/\delta\).
The chain rule gives \(\nabla_z u_j=\nabla_xv_j^\epsilon\). The transformed weak problem is~\eqref{eq:W} with \(\eta=\epsilon/\delta\) and \(\vartheta=x_*/\epsilon\), so uniqueness identifies \(u_j\) with its solution. A change of variables also gives
\[
\begin{aligned}
\frac1{|I_\delta|}\int_{I_\delta}\nabla_xv_j^\epsilon(x)\,dx
&=\int_Q\nabla_z u_j(z)\,dz,\\
\frac1{|I_\delta|}\int_{I_\delta}a^\epsilon(x)\nabla_xv_j^\epsilon(x)\,dx
&=\int_Qa(\vartheta+z/\eta)\nabla_z u_j(z)\,dz.
\end{aligned}
\]
Thus \(F_H=F_\eta\) and \(G_H=G_\eta\), and~\Cref{thm:main} yields~\eqref{eq:PHYS}, which completes the proof.
\end{proof}

\begin{proof}[Proof of \Cref{prop:sharpness}]
With \(m\) and \(r\) as in the statement, define
\[
a(y){:}=\operatorname{diag}\bigl((m+r\sin(2\pi y_1))^{-1},\lambda\bigr).
\]
Since \(m-r>\Lambda^{-1}\) and \(m+r<\lambda^{-1}\), the eigenvalues of \(a(y)\) lie in \([\lambda,\Lambda]\), so \(a\in\mathcal M(\lambda,\Lambda)\). For phase zero, set
\[
D_\eta:=\int_0^1\bigl(m+r\sin(2\pi t/\eta)\bigr)\,dt
=m+\frac{r\eta}{2\pi}\bigl(1-\cos(2\pi/\eta)\bigr)
\]
and
\[
u_1(x)=D_\eta^{-1}\int_0^{x_1}
\bigl(m+r\sin(2\pi t/\eta)\bigr)\,dt,
\qquad
u_2(x)=x_2.
\]
Then \(a(x/\eta)\nabla u_1(x)=D_\eta^{-1}e_1\) and \(a(x/\eta)\nabla u_2(x)=\lambda e_2\). Moreover, \(u_1=0,1\) on \(\Gamma_L,\Gamma_R\), respectively, and \(u_2=0,1\) on \(\Gamma_B,\Gamma_T\), respectively. Hence \(u_j-x_j\in V_j\), and integration of the constant fluxes against \(\nabla\phi\), \(\phi\in V_j\), shows that \(u_1,u_2\) satisfy the two mixed problems. Uniqueness identifies them as the corresponding solutions. Their boundary traces and Fubini's theorem give average gradients \(e_1,e_2\), respectively, and hence \(G_\eta=I\) and \(\mathcal A_\eta=F_\eta=\operatorname{diag}(D_\eta^{-1},\lambda)\).

The zero-mean periodic correctors \(\chi_1(y)=-r\cos(2\pi y_1)/(2\pi m)\) and \(\chi_2(y)=0\) satisfy \(a(y)(e_1+\nabla\chi_1(y))=m^{-1}e_1\) and \(a(y)(e_2+\nabla\chi_2(y))=\lambda e_2\).
Therefore \(\mathcal A=\operatorname{diag}(m^{-1},\lambda)\). Finally, let
\(\eta_N=(N+1/2)^{-1}\), \(N\ge1\). Then
\(\cos(2\pi/\eta_N)=-1\), so \(D_{\eta_N}=m+r\eta_N/\pi>m\). The two diagonal matrices differ only in their first entry; substituting this value of \(D_{\eta_N}\) gives~\eqref{eq:sharpness-exact} and completes the proof.
\end{proof}

\bibliographystyle{amsplain}
\bibliography{reference}
\end{document}